\documentclass[leqno]{article}
\usepackage{verbatim, amsmath, amscd, amssymb}
\usepackage[all]{xy}
\usepackage{xcolor}
\newcommand{\bm}{{\bf m}}
\newcommand{\ba}{{\bf a}}\newcommand{\bb}{{\bf b}}\newcommand{\bc}{{\bf c}}\newcommand{\bn}{{\bf n}}

\newcommand{\rR}{\mathrm{R}}
\newcommand{\rs}{\mathrm{s}}

\newcommand{\rY}{\mathrm{Y}}

\newcommand{\bbF}{\mathbb F}

\newcommand{\bbN}{\mathbb N}

\newcommand{\bbZ}{\mathbb Z}

\newcommand{\der}{\mathrm{der}}

\newcommand{\Dist}{\mathrm{Dist}}

\newcommand{\id}{\mathrm{id}}

\newcommand{\ind}{\mathrm{ind}}

\newcommand{\Lie}{\mathrm{Lie}}

\newcommand{\rsc}{\mathrm{sc}}

\newcommand{\St}{\mathrm{St}}

\newcommand{\lbr}{\begin{bmatrix}}
\newcommand{\rbr}{\end{bmatrix}}
\newcommand{\for}{\bigcirc\kern-2.6ex \because}
\newcommand{\forb}{\bigcirc\kern-2.8ex \because}
\newcommand{\forbb}{\bigcirc\kern-3.0ex \because}
\newcommand{\forbbb}{\bigcirc\kern-3.1ex \because}

\newcommand\pf{\noindent {\it {D\'emonstration.}  }}
\newtheorem{thm}{Th\'eor\`eme.}

\newtheorem{lem}{Lemme.}

\begin{document}
\large
\date{}
\title{
{\bf 
Un scindage du morphisme de Frobenius 
sur l'alg\`{e}bre des distributions
d'un groupe  r\'eductif }
\author{
Michel Gros
\\
Universit\'e de Rennes, Campus de Beaulieu, 
\\
CNRS, IRMAR-UMR 6625, 
\\
35042 Rennes cedex, France
\\
michel.gros@univ-rennes1.fr
\and
K\textsc{aneda} Masaharu\footnote{supported in part by JSPS Grants in Aid for Scientific Research 15K04789 
}
\\
Osaka City University
\\
Department of Mathematics
\\
Osaka 558-8585,
Japan
\\
kaneda@sci.osaka-cu.ac.jp
}
}
\maketitle

\begin{abstract}

Pour un groupe alg\'ebrique semi-simple simplement connexe sur un corps alg\'ebriquement clos de caract\'eristique positive, 
nous avons pr\'ec\'edemment construit un scindage de l'endomorphisme de Frobenius
sur son alg\`ebre des distributions. Nous g\'en\'eralisons  la construction au cas de  
des groupes r\'eductifs   et en d\' egageons  les corollaires correspondants.

For a simply connected semisimple algebraic group 
over an algebraically closed field of positive characteristic
we have already constructed a splitting of the Frobenius endomorphism on
its algebra of distributions.
We generalize the construction to the case of  general reductive groups and derive the corresponding  corollaries.

\end{abstract}

Soient $\Bbbk$ un corps alg\'ebriquement clos
de caract\'eristique  $p>0$,
$G$ un $\Bbbk$-groupe alg\'ebrique r\'eductif  
suppos\'e scind\'e sur $\bbZ$,
$G_\bbZ$ une 
$\bbZ$-forme de Chevalley scind\'ee  de $G$
dont le  rang sera not\'e $\ell$
et
$T_\bbZ$ un tore maximal scind\'e de  
$G_\bbZ$.
Soit
$\Dist(F)$
l'endomorphisme de Frobenius,  d\'efini sur 
$\bbF_p$, de l'alg\`{e}bre des distributions $\Dist(G)$ de $G$.
\'Ecrivons
$\Dist(T_\bbZ)=
\bbZ[\binom{H_1}{n_1},\dots, \binom{H_\ell}{n_\ell}|
n_1,\dots,n_\ell\in\bbN]$
\cite[II.1.12]{J}
et soit
$\mu_0=\prod_{i=1}^\ell\binom{H_i-1}{p-1} \in\Dist(T) :=\Dist(T_\bbZ)\otimes_\bbZ\Bbbk$.
Si $T_1$ d\`{e}signe le noyau de Frobenius de   $T$,
$\mu_0$ est l'unique mesure invariante de  $\Dist(T_1)$ qui soit idempotente
\cite[4.7.1]{GK15}.
Nous allons montrer qu'il existe un homomorphisme de 
 $\Bbbk$-alg\`ebres
$\phi:\Dist(G)\to
\mu_0\Dist(G)\mu_0$
tel que
$\Dist(F)\circ\phi=\id_{\Dist(G)}$.
L'existence d'un tel homomorphisme pour  $G$ semi-simple simplement connexe a \'et\'e \'etablie dans 
\cite[Thm. 1.3]{GK11}.

Soient $R$ le syst\`eme  de racines pour  $G_\bbZ$ relativement \`a $T_\bbZ$
et
$U_{\bbZ,\alpha}$ le sous-groupe radiciel de $G_\bbZ$ correspondant \`a $\alpha\in R$.
Notons  $R^+$ une partie positive de
$R$ et 
$R^\rs$ le sous-ensemble des racines simples de $R^+$.
Soit
$E_\alpha^{(n)}$ 
(resp. $F_\alpha^{(n)}$),
$n\in\bbN$,
la base 
$\bbZ$-lin\'{e}aire standard de
$\Dist(U_{\bbZ,\alpha})$
(resp.
$\Dist(U_{\bbZ,-\alpha})$), $\alpha\in R^\rs$.
Pour
$\bm=(m_{i})\in
\bbN^{\ell}$,
posons
$\binom{H}{\bm}=\prod_{i}\binom{H_i}{m_i}$.
Alors, $\Dist(G_{\bbZ})$ est engendr\'e par les  
$E_\alpha^{(n)}$, $F_\alpha^{(n)}$,
$\alpha\in R^\rs$,
$n\in\bbN$,
$\binom{H}{\bm}$,
$\bm\in
\bbN^{\ell}$.
Nous abr\'evierons
$E_\alpha^{(n)}\otimes1,F_\alpha^{(n)}\otimes1,\binom{H}{\bm}\otimes1\in
\Dist(G_\bbZ)\otimes_\bbZ\Bbbk=\Dist(G)$
respectivement en
$E_\alpha^{(n)},F_\alpha^{(n)},\binom{H}{\bm}$.

Le r\'esultat principal de cet article est, plus pr\'ecis\'ement,  le   
\begin{thm}
Il existe un homomorphisme de  $\Bbbk$-alg\`{e}bres 
$\phi:\Dist(G)\to\mu_0\Dist(G)\mu_0$
tel que 
\[
F_\alpha^{(n)}\mapsto
F_\alpha^{(pn)}\mu_0,
E_\alpha^{(n)}\mapsto
E_\alpha^{(pn)}\mu_0,
\binom{H}{\bb}\mapsto 
\binom{H}{p\bb}\mu_0
\quad
\forall\alpha\in R^\rs,n\in\bbN, \bb\in\bbN^\ell,
\]
et qu'on ait
$\Dist(F)\circ\phi=\id_{\Dist(G)}$.

\end{thm}

Comme dans le cas semi-simple simplement connexe, nous en d\'eduisons ((2.3), Thm.)    la pr\'eservation par   contraction par Frobenius des $G$-modules de l'existence de bonnes filtrations et du caract\`ere injectif ou de basculement pour ceux-ci. Par rapport \`a \cite{GK},   l'obtention   de ces propri\'et\'es   requiert quelques am\'enagements (cf. (1.2), Rem)   m\'eritant d'\^etre   explicit\'es mais aussi un nouvel ingr\'edient de nature diff\'erente. Le lecteur perspicace remarquera en effet  qu'il faut savoir pallier l'absence, dans le pr\'esent cadre g\'en\'eral,  de module de Steinberg. Pour ce faire, nous passerons \`a   une extension centrale ad\'equate  de $G$ dont nous d\'emontrons l'existence dans (2.1), Lem.

\setcounter{equation}{0}\begin{center}$1^\circ$
{Preuve}\end{center}

Pour tout  $\alpha\in R$,
posons
$U_\alpha=U_{\bbZ,\alpha}\otimes_\bbZ\Bbbk$,
$U=\prod_{\alpha\in R^+}U_{-\alpha}$ et 
$U^+=\prod_{\alpha\in R^+}U_{\alpha}$.
Rappelons 
(\cite[Satz 7, p. 22]{J73}) que
$\Dist(U^+)$ 
(resp. $\Dist(U)$)
est engendr\'e par les
$E_\alpha^{(n)}$
(resp.
$F_\alpha^{(n)}$), $\alpha\in R^s$, $n\in\bbN$.
Rappelons \'egalement que la mesure $\mu_0$ commute avec tous les 
\'el\'ements de la forme
$E_\alpha^{(np)}$ et $F_\alpha^{(np)}$,
$\alpha\in R^\rs$, $n\in\bbN$.
Gr\^ace \`a
\cite[1.2]{GK11}, on sait d\'ej\`a que les homomorphismes de $\Bbbk$-alg\`ebres
d\'esir\'es 
$\phi^\pm:\Dist(U^\pm)\to\mu_0\Dist(G)\mu_0$
et
$\phi^0:\Dist(T)\to\mu_0\Dist(G)\mu_0$ existent.
Il faut donc d\'emontrer que  $\phi^\pm$ et  $\phi^0$ s'\'etendent \`a  $\Dist(G)$ tout entier.

\setcounter{equation}{0}
\noindent
(1.1)
Montrons tout d'abord le
\begin{lem}\label{}
Soient
$i\in[1,\ell]$, $a,c\in\bbZ$ et $b\in\bbN$.
Dans $\Dist(T)$, on a les \'egalit\'es suivantes :

{\rm{(i)}} 
$\binom{aH_i}{b}\mu_0=0$
sauf si $p | b$.

{\rm{(ii)}} 
$\phi^0\binom{aH_i}{b}=\binom{aH_i}{pb}\mu_0$.

{\rm{(iii)}} 
$\phi^0\binom{H_i+c}{b}=\binom{H_i+pc}{pb}$.

\end{lem}

\pf
Pour all\'eger, on \'ecrira  dans la suite $H$ ‡ la place de  $H_i$.

(i)
Supposons que $p \nmid b$ et soit
$b=pb'+b''$ la divsion euclidienne de $b$ par $p$ avec donc
$b''\in]0,p[$ .
Remarquons tout d'abord que, $T_1$ d\'esignant comme dans l'introduction le   noyau de Frobenius de $T$, on a
\[
\binom{aH}{b''}=
\frac{aH(aH-1)\dots(aH-b''+1)}{b''!}\in\Dist(T_1),
\]
et donc  
\begin{align}
\binom{aH}{b''}\mu_0
&=
\binom{0}{b''}\mu_0
\quad\text{car $\mu_0$ est une mesure invariante de $\Dist(T_1)$}
\\
&\notag=
0.
\end{align}

Si $b'>0$, alors
\[
\binom{aH}{pb'}\binom{aH}{b''}
=
\sum_{i=0}^{b''}
\binom{pb'+b''-i}{b''}\binom{b''}{i}\binom{aH}{b-i}
=
\binom{pb'+b''}{b''}\binom{aH}{b}
=
\binom{aH}{b}.
\]
D'o\`u l'assertion par (1).

(ii)
On peut supposer $a\ne0$.
Si $a=1$ l'assertion est claire par d\'efinition. Proc\'edons alors par r\'ecurrence.
Si $a>1$,
\begin{align*}
\phi^0\binom{aH}{b}
&=
\phi^0(\sum_{i=0}^b\binom{(a-1)H}{b-i}\binom{H}{i})
\\
&=
\sum_{i=0}^b\phi^0(\binom{(a-1)H}{b-i}\phi^0\binom{H}{i}
\quad\text{car $\phi^0$ est multiplicative sur  $\Dist(T)$
}
\\
&=
\sum_{i=0}^b\binom{(a-1)H}{p(b-i)}\binom{H}{pi}\mu_0
\quad\text{par hypoth\`ese de r\'ecurrence}
\\
&=
\sum_{i=0}^{pb}\binom{(a-1)H}{pb-j)}\binom{H}{j}\mu_0
\quad\text{par (i)}
\\
&=
\binom{aH}{pb}\mu_0.
\end{align*}

Maintenant
\begin{align*}
\phi^0\binom{-H}{b}
&=
\phi^0((-1)^b\binom{H+b-1}{b})
=
(-1)^b\phi^0(\sum_{i=0}^b\binom{b-1}{b-i}\binom{H}{i})
\\
&=
(-1)^b
\sum_{i=0}^b\binom{b-1}{b-i}\binom{H}{pi}\mu_0 ;
\end{align*} 
mais d'autre part
\begin{align*}
\binom{-H}{pb}\mu_0
&=
(-1)^{pb}\binom{H+pb-1}{pb}
\mu_0
=
(-1)^{pb}
\sum_{i=0}^{pb}\binom{pb-1}{pb-i}\binom{H}{i}\mu_0
\\
&=
(-1)^{b}
\sum_{i=0}^{b}\binom{pb-1}{pb-pi}\binom{H}{pi}\mu_0
\quad\text{par (i) de nouveau}
\\
&=
(-1)^{b}
\sum_{i=0}^{b}\binom{p(b-1)+p-1}{p(b-i)}\binom{H}{pi}\mu_0
=
(-1)^{b}
\sum_{i=0}^{b}\binom{b-1}{b-i}\binom{H}{pi}\mu_0,
\end{align*}
et donc
$\phi^0\binom{-H}{b}=\binom{-H}{pb}\mu_0$.
On a alors
\begin{align*}
\phi^0\binom{-(a+1)H}{b}
&=
\phi^0(\sum_{i=0}^b
\binom{-aH}{b-i}\binom{-H}{i})
\\
&=
\sum_{i=0}^b
\binom{-aH}{p(b-i)}\binom{-H}{pi}\mu_0
\quad\text{par hypoth\`ese de r\'ecurrence}
\\
&=
\sum_{j=0}^{pb}
\binom{-aH}{pb-j}\binom{-H}{j}\mu_0
\quad\text{par (i)}
\\
&=
\binom{-aH-H}{pb}\mu_0
=
\binom{-(a+1)H}{pb}\mu_0.
\end{align*}

(iii)
On a
\begin{align*}
\phi^0\binom{H+c}{b}
&=
\phi^0(\sum_{i=0}^b\binom{c}{b-i}\binom{H}{i})
=
\sum_{i=0}^b\binom{c}{b-i}\binom{H}{pi}\mu_0
=
\sum_{i=0}^b\binom{pc}{pb-pi}\binom{H}{pi}\mu_0
\\
&=
\sum_{j=0}^{pb}\binom{pc}{pb-j}\binom{H}{j}\mu_0
\quad\text{par (i)}
\\
&=
\binom{H+pc}{pb}\mu_0.
\end{align*}

\setcounter{equation}{0}
\noindent
(1.2)
Posons maintenant
$B=UT$ and $B^+=U^+T$.
Alors
$\Dist(B)$
(resp. $\Dist(B^+)$) est engendr\'e par les $F_\alpha^{(n)}$
(resp. $E_\alpha^{(n)}$), $\alpha\in R^\rs$, et les
$\binom{H}{\bb}$, $\bb\in\bbN^\ell$.
Prouvons maintenant le  

\begin{lem}

Les homomorphismes $\phi^\pm$ et $\phi^0$
s'\'etendent respectivement en des homorphismes de  $\Bbbk$-alg\`ebres
$\phi^{\geq0}:\Dist(B^+)\to\mu_0\Dist(G)\mu_0$
et
$\phi^{\leq0}:\Dist(B)\to\mu_0\Dist(G)\mu_0$.

\end{lem}

\pf
Pour v\'erifier la premi\`ere assertion, il suffit de v\'erifier que les relations  \cite[Lem. 26.3.D]{Hum}
\begin{align}
F_\alpha^{(a)}\binom{H_i+c}{b}=
\binom{H_i-a\alpha(H_i)+c}{b}
F_\alpha^{(a)},
\quad\alpha\in R^\rs,
i\in[1,\ell], a,b\in\bbN, c\in\bbZ,
\end{align}
sont pr\'eserv\'ees par application de  $\phi$, i.e.,
\[
F_\alpha^{(pa)}\phi^0\binom{H_i+c}{b}=
\{\phi^0\binom{H_i-a\alpha(H_i)+c}{b}\}
F_\alpha^{(pa)}.
\]
Le membre de gauche vaut
\begin{align*}
&
F^{(pa)}\binom{H_i+pc}{pb}\mu_0
\quad\text{par (1.1.iii)}
\\
&=
\binom{H_i-pa\alpha(H_i)+pc}{pb}F^{(pa)}\mu_0
\quad\text{par (1)},
\end{align*}
qui est bien \'egal au membre de droite
gr\^ace de nouveau \`a  1.1 (iii)
plus the fact that $\mu_0$ commutes with $F_\alpha^{(pa)}$.

Il s'ensuit que  $\phi
^{\leq0}$ est multiplicative ; de m\^eme pour $\phi
^{\geq0}$.

{\it{Remarque.}} Le lecteur remarquera ici une diff\'erence avec  le cas semi-simple simplement connexe \cite[Cor. 1.2]{GK11} : il est n\'ecessaire   de mutiplier par $\mu_0$ d\`es cette \'etape de la construction.

\setcounter{equation}{0}
\noindent
(1.3)
Notons $\phi$ le prolongement obtenu  par  $\Bbbk$-lin\'earit\'e de   $\phi^\pm$ et $\phi^0$ \`a  $\Dist(G)=\Dist(U)\otimes\Dist(T)\otimes\Dist(U^+)$ tout entier. On veut maintenant d\'emontrer que ce prolongement est multiplicatif. Posant  
$H_\alpha=[E_\alpha,F_\alpha]$ pour tout $\alpha\in R^\rs$, il suffit de voir que les 
relations 
\cite[Lem. 26.2]{Hum},
\begin{align}
E_\alpha^{(a)}F_\alpha^{(b)}
&=
\sum_{r=0}^{\min\{a,b\}}
F_\alpha^{(b-r)}\binom{H_\alpha+2r-a-b}{r}
E_\alpha^{(a-r)},
\quad
a,b\in\bbN, 
\end{align}
sont pr\'eserv\'ees par 
$\phi$, i.e.  que
\begin{align}
E_\alpha^{(pa)}F_\alpha^{(pb)}
\mu_0
&=
\sum_{r=0}^{\min\{a,b\}}
F_\alpha^{(pb-pr)}\{
\phi\binom{H_\alpha+2r-a-b}{r}
\}
E_\alpha^{(pa-pr)}.
\end{align}

Le membre de gauche de (2) est \'egal \`a
\begin{align*}
&
\sum_{r=0}^{\min\{pa,pb\}}
F_\alpha^{(pb-r)}\binom{H_\alpha+2r-pa-pb}{r}
E_\alpha^{(pa-r)}
\mu_0
\quad\text{par (1)}
\\
&=
\sum_{r=0}^{\min\{pa,pb\}}
F_\alpha^{(pb-r)}
E_\alpha^{(pa-r)}
\binom{H_\alpha+2(pa-r)+2r-pa-pb}{r}
\mu_0
\quad\text{par \cite[Lem. 26.3.D]{Hum}}
\\
&=
\sum_{r=0}^{\min\{pa,pb\}}
F_\alpha^{(pb-r)}
E_\alpha^{(pa-r)}
\binom{H_\alpha+pa-pb}{r}
\mu_0.
\end{align*}
\'Ecrivons maintenant
$H_\alpha=
\sum_{i=1}^\ell
c_iH_i$ avec
$c_i\in\bbZ$.
Alors :
\begin{align}
&
\binom
{H_\alpha+pa-pb}{r}
\mu_0
=
\binom{\sum_{i=1}^\ell
c_iH_i+pa-pb}{r}
\mu_0
\\
\notag&=
\sum_{j=0}^r
\binom{\sum_{i=1}^{\ell-1}
c_iH_i+pa-pb}{r-j}
\binom{c_\ell
H_\ell}{j}
\mu_0
\\
\notag&=
\underset{\substack{j\in\bbN
\\
pj\leq r}}{\sum}
\binom{\sum_{i=1}^{\ell-1}
c_iH_i+pa-pb}{r-pj}
\binom{c_\ell
H_\ell}{pj}
\mu_0
\quad\text{par (1.1.i)}
\\
\notag&=
\dots
\\
\notag&=
\underset{\substack{j_\ell\in\bbN
\\
pj_\ell\leq r}}{\sum}
\dots
\underset{\substack{j_1\in\bbN
\\
pj_1\leq r-pj_\ell-\dots-pj_2}}{\sum}
\binom{
pa-pb}{r-p(j_\ell+\dots+j_2+j_1)}
\binom{c_1
H_1}{pj_1}
\dots
\binom{c_\ell
H_\ell}{pj_\ell}
\mu_0,
\end{align}
qui s'annule sauf si   $p|r$.
Il s'ensuit que le membre de gauche de  
(2) est \'egal \`a
\begin{multline}
\sum_{r=0}^{\min\{a,b\}}
F_\alpha^{(pb-pr)}
E_\alpha^{(pa-pr)}
\\
\sum_{j_\ell=0}^{r}
\dots
\sum_{j_1=0}^{
r-j_\ell-\dots-j_2}
\binom{
a-b}{r-(j_\ell+\dots+j_2+j_1)}
\binom{c_1
H_1}{pj_1}
\dots
\binom{c_\ell
H_\ell}{pj_\ell}
\mu_0.
\end{multline}

D'autre part, le membre de droite de  (2) est \'egal \`a
\begin{align}
\sum_{r=0}^{\min\{a,b\}}
&
F_\alpha^{(pb-pr)}
\phi(\binom{H_\alpha+2r-a-b}{r}
E_\alpha^{(a-r)})
\quad\text{par
(1.2)}
\\
\notag&=
\sum_{r=0}^{\min\{a,b\}}
F_\alpha^{(pb-pr)}
\phi(
E_\alpha^{(a-r)}\binom{H_\alpha+2a-2r+2r-a-b}{r}
)
\\
\notag&=
\sum_{r=0}^{\min\{a,b\}}
F_\alpha^{(pb-pr)}
E_\alpha^{(pa-pr)}
\phi
\binom{H_\alpha+a-b}{r}
\end{align}
avec
\begin{align*}
\phi
\binom{H_\alpha+a-b}{r}
&=
\phi
\binom{\sum_{i=1}^\ell
c_iH_i+a-b}{r}
=\phi
(
\sum_{j=0}^r\binom{
\sum_{i=1}^{\ell-1}
c_iH_i+a-b}{r-j}
\binom{c_\ell
H_\ell}{j}
)
\\
&=
\sum_{j=0}^r
\phi
\binom{
\sum_{i=1}^{\ell-1}
c_iH_i+a-b}{r-j}
\phi\binom{c_\ell
H_\ell}{j}
\\
&=
\sum_{j=0}^r
\phi
\binom{
\sum_{i=1}^{\ell-1}
c_iH_i+a-b}{r-j}
\binom{c_\ell
H_\ell}{pj}
\mu_0
\quad\text{par (1.1.ii)}.
\end{align*}
Finalement, par applications r\'ep\'et\'ees de 1.1(ii), il s'ensuit que 
(4) et (3) sont \'egaux, ce qui termine la preuve du th\'eor\`eme.

\setcounter{equation}{0}\begin{center}$2^\circ$
{Contraction}\end{center}

\setcounter{equation}{0}
\noindent
(2.1)
\'Etant donn\'e un  $G$-module
$M$,
sa {\emph{contraction par Frobenius}} 
$M^\phi=\mu_0M$
admet une structure canonique de 
$(\Dist(G),T)$-module,
et donc de $G$-module
\cite[II.1.20]{J}.
Lorsque $G$ poss\`ede un module de Steinberg
$\St$, Donkin a observ\'e que 
$M^\phi$ est isomorphe \`{a} $((\St\otimes\St\otimes M)^{G_1})^{[-1]}$
comme $G$-modules,
avec
$(\St\otimes\St\otimes M)^{G_1}$ les  $G_1$-invariants de  $\St\otimes\St\otimes M$ (ce qui implique que l'action de $G$ sur $(\St\otimes\St\otimes M)^{G_1}$ se factorise par
$G/G_1$, i.e., par le morphisme de  Frobenius  
$G\to G$), et  
$?^{[-1]}$ la ``d\'etorsion'' via l'isomorphisme
$G/G_1\simeq G$, ce qui est pr\'ecis\'ement l'action via   $\phi$ (voir \cite{GK}  pour plus de d\'etails).
Les r\'esultats  principaux de \cite{GK} 
valent donc pour tout groupe alg\'ebrique r\'eductif 
$G$ d\`es qu'on dispose, sous nos hypoth\`eses, de  modules de Steinberg.

L'existence de tels modules n'\'etant pas garantie pour les $G$ comme dans l'introduction, nous allons donc tout d'abord
relever $G$ en un groupe r\'eductif $\hat G$ qui sera une extension centrale de $G$,  admettra un ``module de Steinberg''
\cite[II. 3.18,  Rmk]{J}  et sera aussi  tel que l'idempotent   
$\mu'_0$ attach\'e \`a  $\hat G$ se projette sur celui,
$\mu_0$, attach\'e \`a $G$.
De mani\`ere g\'en\'erale, pour un \'epimorphisme $G' \rightarrow G$ et $T'$  un tore maximal de $G'$
relevant $T$, la  condition pr\'ec\'edente est garantie si   
$\Lie(T'_\bbZ)$ s'envoie surjectivement sur
$\Lie(T_\bbZ)$.
Si
$\rY(T'')$  d\'esigne le groupe de co-caract\`eres d'un tore  $T''$,
cette derni\`ere condition  est en fait \'equivalente \`a la surjectivit\'e de l'application  $\rY(\hat T)\to\rY(T)$.
La construction suivante, dite d'une $z$-extension de $G$,  est un cas particulier de 
\cite[Prop. 3.1]{MS}.
Nous sommes reconnaissants \`a  ABE Noriyuki de nous avoir signal\'e cette r\'ef\'erence.

\begin{lem}
Il existe une extension centrale
$\hat G$ de $G$ avec 
$\hat G$ reductif scind\'e sur $\bbZ$ et de sous-groupe d\'eriv\'e 
$\hat G^\der$ semi-simple simplement connexe  telle que,  $\hat T$ d\'esignant un tore maximal de 
$\hat G$ scind\'e sur $\bbZ$, l'application
$\rY(\hat T)\to \rY(T)$ soit surjective.

\end{lem}

\pf
Soient
$Z=\cap_{\alpha\in R}\ker\alpha$ le centre de $G$,
$
G^\der$ le sous-groupe d\'eriv\'e de $G$ et 
$T^{\der}=\underset{\substack{
\lambda\in\Lambda
\\
\langle\lambda,\alpha^\vee\rangle=0\, \forall\alpha\in R
}}{\cap}\ker\lambda$ 
un tore maximal  scind\'e sur $\bbZ$
\cite[II.1.18]{J}.
Alors, 
$Z^\der=T^\der\cap Z$ est le centre de $G^\der$.
Soient maintenant $\pi:G^\rsc\to G^\der$ un rev\^etement  simplement connexe de $G^\der$ de tore maximal $T^\rsc=\pi^{-1}(T^\der)$
scind\'e sur  $\bbZ$ \cite[II.1.17]{J} et  
$Z^\rsc=\pi^{-1}(Z^\der)=\cap_{\alpha\in R}\ker(\alpha|_{T^\der}\circ\pi)$
le centre de  $G^\rsc$.

Consid\'erons le morphisme 
$\pi':G^\rsc\times Z\to G$ ;
$(g,z)\mapsto\pi(g)z$.
Comme $\pi'$ induit un isomorphisme
$(T^\rsc\times Z)/\ker\pi'\simeq T$,
on a un isomorphisme $(G^\rsc\times Z)/\ker\pi'\simeq G$
avec $
\ker\pi'
=
\{(g,\pi(g)^{-1})|\pi(g)\in T^\der\cap Z=Z^\der\}
=
\{(g,\pi(g)^{-1})|g\in Z^\rsc\}$.
Si l'on introduit $G^\rsc\times^{Z^\rsc }Z :=(G^\rsc\times Z)/{Z^\rsc }$
pour l'action
$(g,z)z'=(gz',\pi(z')^{-1}z)$,
on obtient donc des isomorphismes compatibles
\begin{align}
\xymatrix{
G^\rsc\times^{Z^\rsc }Z
\ar[rr]^-\sim
&&
G
\\
T^\rsc\times^{Z^\rsc }Z
\ar@{.>}[rr]_-\sim
\ar@{^(->}[]!<0ex,-2ex>;[u]
&&
T.
\ar@{^(->}[]!<0ex,2ex>;[u]
}
\end{align}

Soit maintenant  $\hat Z=T^\rsc\times Z$ et posons
$\hat G :=G^\rsc\times^{Z^{\rsc}}\hat Z=
(G^\rsc\times(T^\rsc\times Z))/Z^{\rsc}$ pour l'action
$(g,(t,z'))z=(gz, (z^{-1}t, \pi(z^{-1})z'))$,
$g\in G^\rsc,t\in T^\rsc,z'\in Z, z\in Z^\rsc$.
On a une immersion ferm\'ee
$\hat Z\hookrightarrow\hat G$ via $(t,z)\mapsto[1,(t,z)]$
\cite[I.5.6]{J}.
D\'efinissons  le morphisme
$\eta:\hat G\to G$ par
$[g,(t.z')]\mapsto\pi(g)z'$;
pour tout $ z\in Z^\rsc$,
$[gz,(z^{-1}t.\pi(z^{-1})z')]\mapsto\pi(gz)\pi(z^{-1})z'=\pi(g)z'$.
Utilisant l'inclusion $T^\rsc\hookrightarrow\hat Z$
donn\'ee par
$t\mapsto(t,1)$, on voit que 
$\eta$ se factorise \`a travers 
$\hat G/T^\rsc\simeq G^\rsc\times^{Z^\rsc}Z$
via
$[g,(t,z)]T^\rsc\mapsto[g,z]$ d'inverse
$[g,z]\mapsto[g,(1,z)]T^\rsc$.
On obtient donc \`a partir de (1)
un diagramme commutatif de suites exactes de groupes alg\'ebriques  
\begin{equation}
\xymatrix{
1
\ar[r]
&
T^{\rsc}
\ar[r]
& 
\hat G
\ar[r]
&
G
\ar[r]
&
1
\\
1
\ar[r]
&
T^{\rsc}
\ar[r]
\ar@{=}[u]
& 
T^\rsc\times^{Z^\rsc}\hat Z
\ar[r]
\ar@{^(->}[]!<0ex,-2ex>;[u]
&
T
\ar[r]
\ar@{^(->}[]!<0ex,2ex>;[u]
&
1.
}
\end{equation}

En particulier, $\hat G$ est connexe.
Finalement, $\hat G$ est un groupe r\'eductif de sous-groupe d\'eriv\'e isomorphe \`a $G^\rsc$ via le plongement $G^\rsc\hookrightarrow \hat G$
d\'efini par
$g\mapsto[g,(1,1)]$.

De plus,
la suite exacte du bas de  (2) se scinde :
on a en effet une identification $T^\rsc\times^{Z^\rsc}\hat Z\simeq
 T^\rsc\times(T^\rsc\times^{Z^\rsc}Z)$
via
$[t_1,(t_2,z)]\mapsto(t_1t_2,[t_1,z])$
d'inverse
$(t_1,[t_2,z])\mapsto[t_2,(t_1t_2^{-1},z)]$.
Ainsi, 
$T^\rsc\times^{Z^\rsc}\hat Z$ est un tore maximal de
$\hat G$ scind\'e sur $\bbZ$
et l'on a  une surjection 
$\rY(T^\rsc\times^{Z^\rsc}\hat Z)\twoheadrightarrow\rY(T)$  comme cherch\'e.

\setcounter{equation}{0}
\noindent
(2.2)
Conservons les hypoth\`eses et notations de  (2.1).
Soit
$\hat\pi:\hat G\to G$ l'extension centrale du lemme pr\'ec\'edent et   posons $K=\ker\hat\pi=T^{\rsc}$.
D\'efinissons  l'idempotent
$\hat\mu_0\in \Dist(\hat T_1)$
de la m\^eme mani\`ere que pour $T$, et soit $\hat\phi:\Dist(\hat G)\to \hat\mu_0\Dist(\hat G)\hat\mu_0$
le scindage de Frobenius  correspondant $\hat F$ sur $\hat G$.

\begin{lem}
{\rm{(i)}}
On a un diagramme commutatif de  
$\Bbbk$-alg\`ebres
\[
\xymatrix{
\Dist(\hat G)
\ar[d]_-{\Dist(\hat\pi)}
\ar[rr]^-{\hat \phi}
&&
\hat\mu_0\Dist(\hat G)\hat \mu_0
\ar[d]^-{\Dist(\hat\pi)
}
\\
\Dist(G)
\ar[rr]_-\phi
&&
\mu_0\Dist(G)\mu_0
}
\]
avec
$\Dist(\hat\pi)(\hat\mu_0)=\mu_0$.

{\rm{(ii)}}
Pour un un $\hat G$-module 
$M$
on a un isomorphisme de  $G$-modules
\[
(M^K)^{\phi}\simeq
(M^{\hat\phi})^K.
\]

\end{lem}

\pf
(i)
Soit
$(\hat H_i|i\in[1,\ell'])$ 
une $\bbZ$-base de $\Lie(\hat T_\bbZ)$
relevant la
 $\bbZ$-base  
$(H_i|i\in[1,\ell])$
de $\Lie(T_\bbZ)$
correspondant \`a
$Y(\hat T)\twoheadrightarrow
Y(T)$ ; 
\[
\Dist(\hat\pi)(\hat H_i)
=
\begin{cases}
H_i
&\text{si $i\in[1,\ell]$},
\\
0
&\text{sinon}.
\end{cases}
\]
Comme $\binom{-1}{p-1}=1$, on a en particulier,
\begin{align}
\Dist(\hat\pi)(\hat\mu_0)=\mu_0.
\end{align}
Alors, pour tous $i\in[1,\ell']$ et $m\in\bbN$,

\begin{align*}
\mu_0\Dist(\hat\pi)(\hat\phi\binom{\hat H_i}{m})\mu_0
&=
\binom
{\Dist(\hat\pi)
(\hat H_i)}{pm}
\mu_0
=
\begin{cases}
\binom
{H_i}{pm}
\mu_0
&\text{si $i\in[1,\ell]$},
\\
0
&\text{sinon}
\end{cases}
\\
&=
\phi(\Dist(\hat\pi)\binom{\hat H_i}{m}).
\end{align*}
L'assertion s'ensuit.

(ii)
L'ensemble sous-jacent \`a ces deux $G$-modules est
\[
(M^K)^{T_1}=
(M^K)^{\hat T_1}
=
M^{K\hat T_1}
=
(M^{\hat T_1})^K.
\]
D'apr\`es  (i), 
l'isomorphisme annonc\'e est donc juste l'identit\'e.

\setcounter{equation}{0}
\noindent
(2.3)
Rappelons (\cite[I.6.11]{J}) 
que pour tout $\lambda\in\Lambda$, on a 
\[
\rR^\bullet\ind_{\hat B}^{\hat G}(\lambda\circ\hat\pi)
=
\{
\rR^\bullet\ind_{\hat B}^{\hat G}(\lambda\circ\hat\pi)
\}^K
\simeq
\rR^\bullet\ind_{\hat B/K}^{\hat G/K}(\lambda\circ\hat\pi)
\simeq
\rR^\bullet\ind_{B}^{G}(\lambda).
\]
et aussi 
(\cite[I.6.4]{J}) que si  $I$ est un $\hat G$-module injectif,
$I^K$ demeure injectif comme $G$-module.
Comme
$\Bbbk[\hat G]^K=\Bbbk[\hat G/K]=\Bbbk[G]$, on a $\Bbbk[G]^\phi\simeq
(\Bbbk[\hat G]^{\hat\phi})^K$.
Tout $G$-module injectif \'etant un facteur direct  d'une somme directe de copies de   
$\Bbbk[G]$
\cite[I.3.9]{J}, du cas dans lequel  
$G$ admet un module de  Steinberg  \cite{GK}/\cite{A}, on   d\'eduit le 
\begin{thm}
La contraction par Frobenius pour les $G$-modules pr\'eserve
l'existence de bonnes filtrations et envoie modules de basculement sur modules de basculement et modules injectifs sur modules injectifs.
\end{thm}

\end{document}